# A Simple Linear-Time Algorithm for Finding Path-Decompositions of Small Width

Kevin Cattell, Michael J. Dinneen, Michael R. Fellows

Department of Computer Science
University of Victoria
Victoria, British Columbia V8W 3P6, Canada
`mjd,kcattell,mfellows @csr.uvic.ca`

## Abstract

We described a simple algorithm running in linear time for each fixed constant $k$, that either establishes that the pathwidth of a graph $G$ is greater than $k$, or finds a path-decomposition of $G$ of width at most $O(2^k)$. This provides a simple proof of the result by Bodlaender that many families of graphs of bounded pathwidth can be recognized in linear time.

**Classification:** Algorithms and data structures, computational complexity.

## 1 Introduction

The topics of the *pathwidth* and *treewidth* of graphs have proven to be of fundamental interest for two reasons. First of all, they play an important role in the deep results of Robertson and Seymour [RS83, RS86a, RS90, RS91, RS]. Secondly, and more importantly from a practical point of view, bounded pathwidth and treewidth have proven to be general "common denominators" for many natural input restrictions of NP-complete problems. For many important problems, we now know that fixing a natural parameter $k$ implies that the yes-instances have bounded treewidth or pathwidth (for examples see [Bod88b, BFW92, FHW93, FL92, KT92, Moh90]). We also know that many problems can be solved in linear time when the input includes a bounded-width path-decomposition (or tree-decomposition) of the graph (see [Arn85, ALS91, AP89, Bod88a, BPT92, CM93, WHL85, Wi87] and [Bod92] for many further references).

After several rounds of improvement [RS86b, La90, Re92] the best known algorithm for finding tree-decompositions is due to Bodlaender [Bod93]. For each fixed $k$, this algorithm in time $O(2^{k^2}n)$ either determines that the treewidth is greater than $k$, or produces a tree-decomposition of width at most $k$. By first running this algorithm and then applying the algorithm of [BK91, Kl93], a similar result holds for pathwidth. Both of the algorithms involved are quite complicated.

We describe here a very simple algorithm based on "pebbling" the graph using a pool of $O(2^k)$ pebbles, that in linear time (for fixed $k$), either determines that the pathwidth of a graph is more than $k$, or finds a path-decomposition of width at most the number of pebbles actually used. The main advantages of this algorithm over previous results are: (1) the simplicity of the algorithm and (2) the improvement of the hidden constant for a determination that the pathwidth is greater than $k$. The main disadvantage is in the width of the resulting "approximate" decomposition when the width is less than or equal to $k$.

## 2 Preliminaries

All of our discussion concerns finite simple graphs. Some of the graphs have a *boundary* of size $k$, meaning that they have a distinguished set of vertices labeled $1, 2, \ldots, k$. Two boundaried graph can be glued (to form a regular graph) with the $\oplus$ operator, which simply identifies vertices with the same boundary label.

An *(homeomorphic) embedding* of a graph $G_1 = (V_1, E_1)$ in a graph $G_2 = (V_2, E_2)$ is an injection from vertices $V_1$ to $V_2$ with the property that the edges $E_1$ are mapped to disjoint paths of $G_2$. (These disjoint paths in $G_2$ represent possible *subdivisions* of the edges of $G_1$.) The set of homeomorphic embeddings between graphs gives a partial order, called the *topological order*.

A *lower ideal* $\mathcal{J}$ in a partial order $(\mathcal{U}, \geq)$ is a subset of $\mathcal{U}$ such that if $X \in \mathcal{J}$ and $X \geq Y$ then $Y \in \mathcal{J}$. The *obstruction set* for $\mathcal{J}$ is the set of minimal elements of $\mathcal{U} - \mathcal{J}$.

**Definition.** A *path-decomposition* of a graph $G = (V, E)$ is a sequence $X_1, X_2, \ldots, X_r$ of subsets of $V$ that satisfy the following three conditions:

1. $\bigcup_{1 \leq i \leq r} X_i = V$,

2. for every edge $(u, v) \in E$, there exists an $X_i$ such that $u \in X_i$ and $v \in X_i$, and

3. for $1 \leq i < j < k \leq r$, $X_i \cap X_k \subseteq X_j$.

The *pathwidth of a path-decomposition* $X_1, X_2, \ldots, X_r$ is $\max_{1 \leq i \leq r} |X_i| - 1$. The *pathwidth of a graph* $G$ is the minimum pathwidth over all path-decompositions of $G$. Determining pathwidth is equivalent to several VLSI layout problems such as *gate matrix layout* and *vertex separation* [Moh90, EST87].

It is easy to see that the family of graphs of pathwidth at most $t$ is a lower ideal in the topological (and minor) order. It is also known that those graphs with order $n$ have at most $nt - (t^2 + t)/2$ edges.

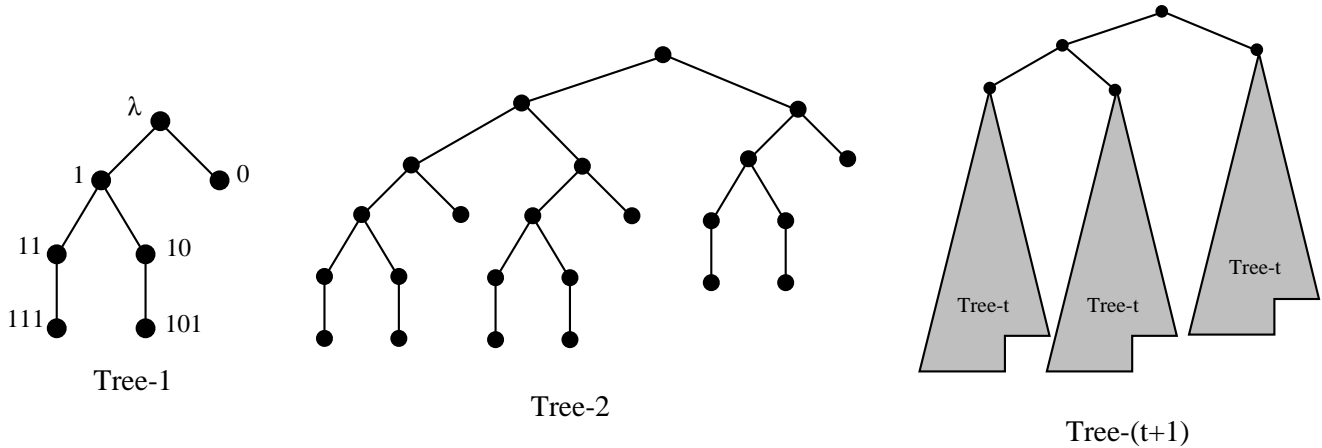

Figure 1: Embedding pathwidth $t$ tree obstructions Tree-$t$ in binary trees.

Let $B_h$ denote the complete binary tree of height $h$ and order $2^h - 1$. Let $h(t)$ be the least value of $h$ such that $B_{h(t)}$ has pathwidth greater than $t$, and let $f(t)$ be the number of vertices of $B_{h(t)}$. To get a bound for $f(t)$, $B_{h(t)}$ needs to contain at least one obstruction of pathwidth $t$. In [EST87] it is shown that all topological tree obstructions of pathwidth $t$ can be recursively generated by the following rules.

1. The single edge tree $K_2$ is the only obstruction of pathwidth 0.

2. If $T_1, T_2$ and $T_3$ are any 3 tree obstructions for pathwidth $t$ then the tree $T$ consisting of a new degree 3 vertex attached to any vertex of $T_1, T_2$ and $T_3$ is a tree obstruction for pathwidth $t + 1$.

From this characterization we see that the orders the tree obstructions of pathwidth $t$ are precisely $(5 \cdot 3^t - 1)/2$, (e.g., orders 2, 7, 22 and 57 for pathwidth $t = 0, 1, 2$ and 3). We can easily embed at least one of the tree obstructions for pathwidth $t$, as shown in Figure 2, in the complete binary tree of height $2t + 2$. Thus, the complete binary tree of order $f(t) = 2^{2t+2} - 1$ has pathwidth greater than $t$.

## 3 Pathwidth Algorithm

Using the $f(t)$ bound given in the previous section, the main result of the paper now follows:
**Theorem 1.** Let $H$ be an arbitrary undirected graph, and let $t$ be a positive integer. One of the following two statements must hold:
(a) The pathwidth of $H$ is at most $f(t) - 1$.
(b) $H$ can be factored: $H = A \oplus B$, where $A$ and $B$ are boundaried graphs with boundary size $f(t)$, the pathwidth of $A$ is greater than $t$ and less than $f(t)$.

**Proof.** We describe a algorithm that terminates either with a path-decomposition of $H$ of width at most $f(t) - 1$, or with a path-decomposition of a suitable factor $A$ with the last vertex set of the decomposition consisting of the boundary vertices.

If we find an homeomorphic embedding of the *guest tree* $B_{h(t)}$ in the *host graph* $H$ then we know that the pathwidth of $H$ is greater than $t$. During the search for such an embedding, we work with a partial embedding. We refer to the vertices of $B_{h(t)}$ as *tokens*, and call tokens *placed* or *unplaced* according to whether or not they are mapped to vertices of $H$ in the current partial embedding. A vertex $v$ of $H$ is *tokened* if a token maps to $v$. At most one token can be placed on a vertex of $H$ at any given time. We recursively label the tokens by the following standard rules:

1. The root token is labeled by the empty string $\lambda$.

2. The left child token and right child token of a height $h$ parent token $P = b_1 b_2 \cdots b_h$ are labeled $P \cdot 1$ and $P \cdot 0$, respectively.

Let $P[i]$ denote the set of vertices of $H$ that have a token at time step $i$. The sequence $P[0], P[1], \ldots, P[s]$ will describe a path-decomposition either of the entirety of $H$ or of a factor $A$ fulfilling the conditions of Theorem 1. In the case of outcome (b) the boundary of the factor $A$ is indicated by $P[s]$.

The placement algorithm is described as follows. Initially consider that every vertex of $H$ is colored *blue*. In the course of the algorithm a vertex of $H$ has its color changed to *red* when a token is placed on it, and stays red if the token is removed. Only blue vertices can be tokened, and so a vertex can only be tokened once.

    **function** GrowTokenTree
1    **if** root token $\lambda$ is not placed on $H$ **then**
        arbitrarily place $\lambda$ on a blue vertex of $H$
    **endif**
2    **while** there is a vertex $u \in H$ with token $T$ and blue neighbor $v$,
        and token $T$ has an unplaced child $T \cdot b$ **do**
    2.1  place token $T \cdot b$ on $v$
    **endwhile**
3    **return** $\{$tokened vertices of $H\}$

    **program** PathDecompositionOrSmallFatFactor
1    $i \leftarrow 0$
2    $P[i] \leftarrow$ **call** GrowTokenTree
3    **until** $|P[i]| = f(t)$ **or** $H$ has no blue vertices **repeat**
    3.1  pick a token $T$ with untokened children
    3.2  remove $T$ from $H$
    3.3  **if** $T$ had one tokened child **then**

    replace all tokens $T \cdot b \cdot S$ with $T \cdot S$
  **endif**
 3.4 $i \leftarrow i + 1$
 3.5 $P[i] \leftarrow$ **call** GrowTokenTree
**enduntil**
**done**

Before we prove the correctness of the algorithm, we note some properties: (1) the root token will need to be placed (step 1 of the GrowTokenTree) at most once for each component of $H$; (2) the GrowTokenTree function only returns when $B_{h(t)}$ has been embedded in $H$ or all parent tokens of degree less than 2 have no blue neighbors; (3) the algorithm will terminate since during each iteration of step 3.2 a tokened red vertex becomes untokened, and this can happen at most $n$ times, where $n$ denotes the order of the host $H$.

Since tokens are placed only on blue vertices and are removed only from red vertices, it follows that the interpolation property of a path-decomposition is satisfied. Suppose the algorithm terminates at time $s$ with all of the vertices colored red. To see that the sequence of vertex sets $P[0], \ldots, P[s]$ represents a path-decomposition of $H$, it remains only to verify that for each edge $(u, v)$ of $H$ there is a time $i$ with both vertices $u$ and $v$ in $P[i]$. Suppose vertex $u$ is tokened first and untokened before $v$ is tokened. But vertex $u$ can be untokened only if all neighbors, including vertex $v$, are colored red (see step 3.1 and comment (2) above).

Suppose the algorithm terminates with all tokens placed. The argument above establishes that the subgraph $A$ of $H$ induced by the red vertices, with boundary set $P[s]$ has pathwidth at most $f(t)$. To complete the proof we argue that in this case the sequence of token placements establishes that $A$ contains a subdivision of $B_{h(t)}$, and hence must have pathwidth greater than $t$. Since the GrowTokenTree function only attaches pendant tokens to parent tokens we need to only to observe that the operation in step 3.3 subdivides the edge between $T$ and its parent. □

**Corollary 1.** Given a graph $H$ of order $n$ and an integer $t$, there exists a $O(n)$ time algorithm that gives evidence that the pathwidth of $H$ is greater than $t$ or finds a path-decomposition of width at most $O(2^t)$.

**Proof.** We show that program PathDecompositionOrSmallFatFactor runs in linear time. First, if $H$ has more than $t \cdot n$ edges, then the pathwidth of $H$ is greater than $t$. By the proof of Theorem 1, the program terminates with either the embedded binary tree as evidence, or a path-decomposition of width at most $f(t)$.

Note that the guest tree $B_{h(t)}$ has constant order $f(t)$, and so token operations that do not involve scanning $H$ are constant time. In function GrowTokenTree, the only non-constant time operation is the check for blue neighbors in step 2. While scanning the adjacent edges of vertex $u$ any edge to a red vertex can be removed, in constant time. Edge $(u, v)$ is also removed when step 2.1 is executed. Therefore, across all calls to GrowTokenTree, each edge of $H$ needs to be considered at most once, for a total of $O(n)$ steps. In program PathDecompositionOrSmallFatFactor, all steps except for GrowTokenTree are constant time. The total number of iterations through the loop is bounded by $n$, by the termination argument following the program. □

The next result shows that we can improve the pathwidth algorithm by restricting the guest tree. This allows us to use the subdivided tree obstructions given in Figure 2.

**Corollary 2.** Any subtree of the binary tree $B_{h(t)}$ that has pathwidth greater than $t$ may be used in the algorithm for Theorem 1.

**Proof.** The following simple modifications allow the algorithm to operate with a subtree. The subtree is specified by a set of *flagged* tokens in $B_{h(t)}$. At worst, the algorithm can potentially embed all of $B_{h(t)}$.

In step 2 of GrowTokenTree, the algorithm only looks for a flagged untokened child $T \cdot b$ to place, since unflagged tokens need not be placed. The stopping condition in step 3 of PathDecompositionOrSmallFatFactor is changed to "all flagged tokens of $B_{h(t)}$ are placed **or** ...," so that termination occurs as soon as the subtree has been embedded. The relabeling in step 3.3 can place unflagged tokens on vertices of $H$ since all the rooted subtrees of a fixed height are not isomorphic. If that happens, we expand our guest tree with those new tokens. (It is easy to see that the new guest is still a tree.) These tokens can then get relabeled by future edge subdivisions that occur above the token in the host tree. Thus, duplication of token labels will not happen, and the $f(t)$ width bound is preserved. □

**Example.** Using $K_{1,3}$ as the guest tree of pathwidth 2 the following program trace terminates with all vertices colored red (gray) yielding a path-decomposition of width 5.

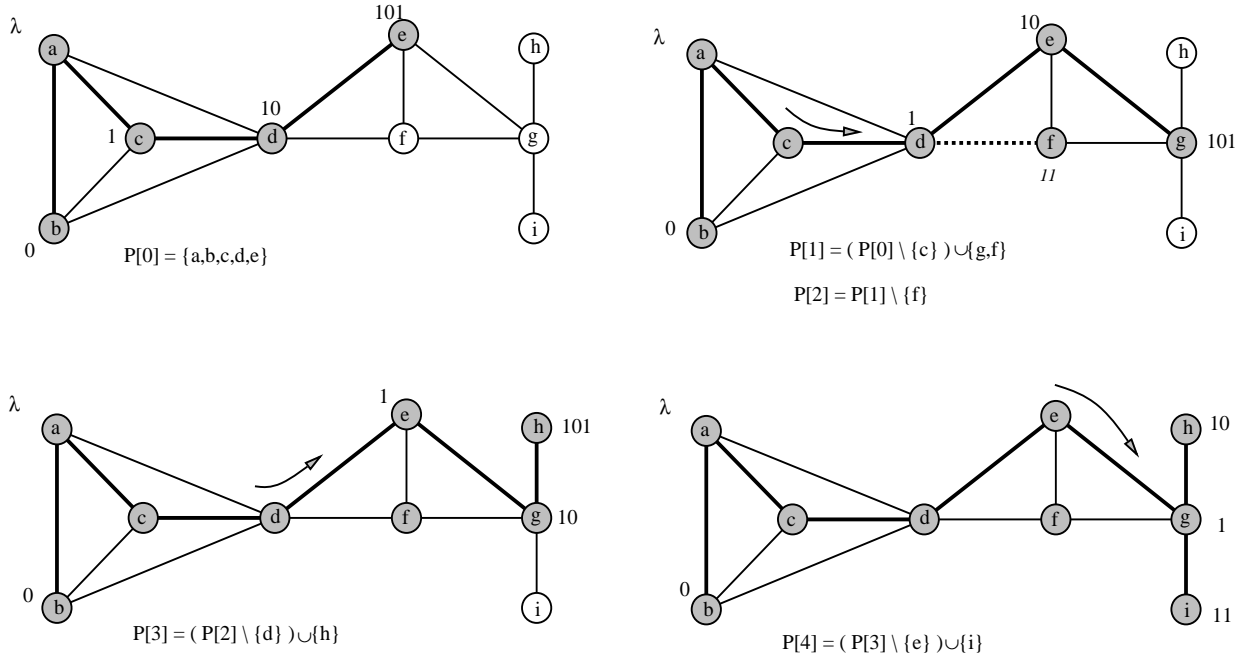

In the proof of Corollary 2 one may wish to not expand the guest tree by flagging new tokens. This can be done and, in fact, is what we would do in practice. Without loss of generality, suppose token $T \cdot 1$ is on vertex $u \in H$ and has children that can not be placed on $H$, and $T \cdot 1$ has one unflagged sibling token $T \cdot 0$ on $v \in H$. If we ignore the flagging of new vertices in the current algorithm, the token $T \cdot 1$ would be removed and the parent $T$ (which has only one legitimate child) would be placed on vertex $u$. What happens to any blue vertices of adjacent to only vertex $u$ (or its unflagged subtree)? The answer is that they are lost and the algorithm would not terminate unless it could embed the guest tree in the remaining portion of $H$. We can fix this problem by checking for unflagged siblings before step 3.3 and to shift the token $T \cdot 1$ from $u$ to $v$. See step $3.3'$ below.

> $3.3'$ **if** $T = P \cdot b$ had an unflagged sibling **then**
>    replace all tokens $P \cdot b \cdot S$ with $P \cdot \text{not}(b) \cdot S$
> **else if** $T$ had one tokened child **then**
>    replace all tokens $T \cdot b \cdot S$ with $T \cdot S$
> **endif**

# 4 Further Directions

In the case that the pathwidth of an input graph $G$ is at most $k$, our algorithm yields a path-decomposition that can have a width exponential in $k$, but that is equal in any case to the maximum number of tokens placed on the graph at any given time, minus 1. It would be interesting to know if this exponential bad behavior is "normal" or whether the algorithm tends to use a smaller number of tokens in practice. Since the pebbling proceeds according to a greedy strategy with much flexibility, there may be placement heuristics that can improve its performance on "typical" instances.